\documentclass{article}
\usepackage[utf8]{inputenc}
\usepackage[english]{babel}
\usepackage{amsmath}
\usepackage{amsfonts}
\usepackage{graphicx}
\usepackage{tikz}
\newtheorem{theorem}{Theorem}
\newtheorem{lemma}{Lemma}

\usepackage{hyperref}
\newcommand{\footremember}[2]{%
	\footnote{#2}
	\newcounter{#1}
	\setcounter{#1}{\value{footnote}}%
}
 
\title{Not Conway's 99-Graph Problem}
\author{%
	Sa'ar Zehavi\footremember{alley}{Technion, Department of Computer Science, saarzehavi@gmail.com}%
	\and Ivo Fagundes David de Oliveira\footremember{trailer}{Technion, Department of Computer Science, ivodavid@gmail.com}%
}

\begin{document}
\maketitle

\section{Abstract}
	Conway's 99-graph problem is the second problem amongst the five 1000\$ 2017 open problems set~\cite{CON}. Four out of the five remain unsolved to this day, including the 99-graph problem. In this paper we quote Conway's definition of the problem and give an alternative interpretation of it, which we humorously name ``not Conway's 99-graph problem". We solve the alternative interpretation completely.

\section{Introduction}
	Conway's 99-graph problem~\cite{CON} asks the following: ``Is there a graph with 99 vertices in which every edge (i.e. pair of joined vertices)	belongs to a unique triangle and every none-edge (pair of unjoined vertices) to a unique quadrilateral?". 
	
	Put formally, is there a simple, undirected graph $G=(V,E)$, such that the following holds?
	\begin{itemize}
		\item ($\bar{P}1$): $|V| = 99$
		\item ($\bar{P}2$): $\forall \{x,y\}\in E$, there exists a unique $z\in V$, such that\\ $\{x,z\},\{y,z\}\in E$.
		\item ($\bar{P}3$): $\forall \{x,y\}\in E^c$, there exists unique pair $z,w\in V$, such that\\ $\{x,z\},\{y,z\},\{x,w\},\{y,w\}\in E$.
	\end{itemize}
	
	It can be proven that a graph satisfying ($\bar{P}2-3$) is, in fact, regular. Whether a regular graph, satisfying ($\bar{P}1-3$) exists, is a known open problem. Such a graph is also called strongly regular. Strongly regular graphs are a parametric family of graphs. A counting argument shows that a strongly regular graph satisfying ($\bar{P}1-3$) has vertex degree of 14.
	
	An alternative interpretation of Conway's 99-graph problem is the following, which we humorously name ``not Conway's 99-graph problem":
	
	Is there a simple, undirected graph $G=(V,E)$, such that the following holds?
		\begin{itemize}
		\item (P1): $|V| = 99$
		\item (P2): $\forall e\in E$, there exists a unique triangle $\Delta\in G$ such that $e\in \Delta$.
		\item (P3): $\forall e\in E^c$, there exists a unique quadrilateral $\Box\in G^c$ such that $e\in \Box$.
	\end{itemize}
	
	We prove the following:
	\begin{theorem}
		\label{theorem:conway}
		Assume $G=(V,E)$, is a simple graph satisfying (P2) and (P3), then $|V| = 5$ or $|V| = 3$. I.e. $G$ is either a triangle or two triangles intersecting in a single vertex.
	\end{theorem}

	The previous theorem, in particular shows that there does not exist a $G$ satisfying (P1-3) simultaneously. In the next section we supply an overview of the proof of our theorem.
\section{Overview}
	For simplicity, due to our use of illustrations, let us name the edges ``blue edges", and the none-edges ``red edges". Assuming $G$ satisfies (P2-3), we know that every blue edge belongs to a unique blue triangle, and every red edge belongs to a unique red quadrilateral. Our main 2 lemmas are the following.
	
	\begin{lemma}
		\label{lemma:one}
		If $G=(V,E)$ satisfies (P2-3), then $G$ has no red triangles.
	\end{lemma}

	\begin{lemma}
		If $G=(V,E)$ satisfies (P2-3), then $G$ has no blue quadrilaterals.
	\end{lemma}

	Assuming $G$ is a graph with $n$ vertices, satisfying (P2-3), the previous lemmas imply that the graph induced by the red edges is triangle free and hence has at most $n^2/4$ edges, and that the graph induced by the blue edges is quadrilateral free, and hence has at most $(1 + \sqrt{2n - 3})n/4$ edges. The red and blue edges, together, form the complete graph, and hence the following inequality must hold: $n^2/4 + (1 + \sqrt{2n - 3})n/4 \ge {n\choose 2}$. The previous implies that $n<9$.
	
	\begin{theorem}
		\label{theorem:main}
		If $G$ satisfies (P2-3), then $|V|<9$.
	\end{theorem}
	Theorem~\ref{theorem:main} implies the inexistence of a graph satisfying (P1-3), concluding with a \emph{negative} answer to the ``not 99-graph problem". Theorem~\ref{theorem:main}, together with extensive search (for graphs with at most 8 vertices) implies Theorem~\ref{theorem:conway}, which states the only possible configurations for graphs satisfying (P2-3) are attained when $n=3$ or $n=5$, hence solving the ``not 99-graph problem" for all possible values of $n$.
	
\section{$G$ has no red triangles - proof of Lemma 1}
	We will prove by case analysis that $G$ satisfying (P2-3) has no red triangles. In the following, there are several graph illustrations. Red and blue edges are denoted by colored lines. Every edge in our graphs is either red or blue. In case an edge does not appear in our illustrations it is either undecided or irrelevant to the argument. Let us assume by contradiction the existence of a red triangle, then there is some red triangle $\Delta(ABC)$:
	
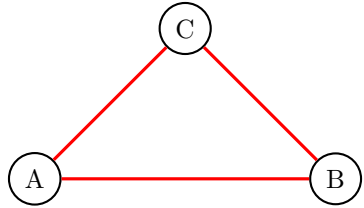
\begin{center}
		\begin{tikzpicture}
			\begin{scope}[every node/.style={circle,thick,draw}]
				\node (A) at (0,0) {A};
				\node (B) at (4,0) {B};
				\node (C) at (2,2) {C};
			\end{scope}		
			\begin{scope}[every edge/.style={draw=red,very thick}]
				\path (A) edge [above] node {} (B);
				\path (B) edge [above] node {} (C);
				\path (C) edge [above] node {} (A);
			\end{scope}
		\end{tikzpicture}
	\end{center}
	As every red edge is a part of a unique red quadrilateral, there is some red quadrilateral containing $AB$. As a quadrilateral cannot contain all edges of a triangle, there are two possible cases:
	\begin{itemize}
		\item Case 1 - Each edge of $\Delta(ABC)$ lies in a unique quadrilateral.
		\item Case 2 - There are two edges of $\Delta(ABC)$ sharing a quadrilateral.
	\end{itemize}

	\begin{lemma}
		Case 1 is impossible.
	\end{lemma}

	\emph{Proof}: Let us assume the first case holds, then there are vertices $D$ and $E$, different than $A$, $B$ and $C$, such that $AB\in\Box(ABED)$:
	
	* Note that quadrilateral $\Box(ABED)$, refers to the cycle of length 4 comprised of the sides $AB,BE,ED,DA$.
	
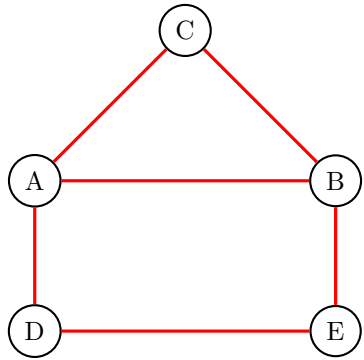
\begin{center}
		\begin{tikzpicture}
			\begin{scope}[every node/.style={circle,thick,draw}]
				\node (A) at (0,0) {A};
				\node (B) at (4,0) {B};
				\node (C) at (2,2) {C};
				\node (D) at (0,-2) {D};
				\node (E) at (4,-2) {E};
			\end{scope}		
			\begin{scope}[every edge/.style={draw=red,very thick}]
				\path (A) edge [above] node {} (B);
				\path (B) edge [above] node {} (C);
				\path (C) edge [above] node {} (A);
				\path (A) edge [above] node {} (D);
				\path (B) edge [above] node {} (E);
				\path (D) edge [above] node {} (E);
			\end{scope}
		\end{tikzpicture}
	\end{center}
	
	The edge $AC$ also lies in a unique quadrilateral, which by our assumption, is edge wise disjoint to $\Box(ABED)$. Therefore, there are vertices $F,G\in V$, such that $AC$ is in the red quadrilateral $\Box(ACGF)$. Is it possible that $\{D,E\}\cap\{F,G\}\neq \emptyset$? Assume that $D$ is shared by the two quadrilaterals. This implies that $\Box(ACGF)$ is either $\Box(ACDF)$ or $\Box(ACGD)$. Clearly, as $\Box(ACGD)$ and $\Box(ABED)$ share $AD$, in contradiction to (P3), it has to be that $D\neq F$. Thus, we collapse to $\Box(ACGF) = \Box(ACDF)$, in which case $CD$ is red, and we have:
	
	
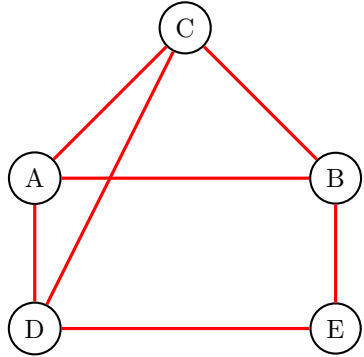
\begin{center}
		\begin{tikzpicture}
			\begin{scope}[every node/.style={circle,thick,draw}]
				\node (A) at (0,0) {A};
				\node (B) at (4,0) {B};
				\node (C) at (2,2) {C};
				\node (D) at (0,-2) {D};
				\node (E) at (4,-2) {E};
			\end{scope}		
			\begin{scope}[every edge/.style={draw=red,very thick}]
				\path (A) edge [above] node {} (B);
				\path (B) edge [above] node {} (C);
				\path (C) edge [above] node {} (A);
				\path (A) edge [above] node {} (D);
				\path (B) edge [above] node {} (E);
				\path (D) edge [above] node {} (E);
				\path (C) edge [above] node {} (D);
			\end{scope}
		\end{tikzpicture}
	\end{center}
	
	This is impossible since now $BE$ is shared by both $\Box(ABED)$ and $\Box(CBED)$, a contradiction to (P3). Similarly, we reach a contradiction by assuming $E\in\{F,G\}$, as this implies either $EA$ or $EC$ are red. Thus, assuming we are in case one implies that $\{F,G\}\cap\{D,E\}=\emptyset$. So we have:
	
	
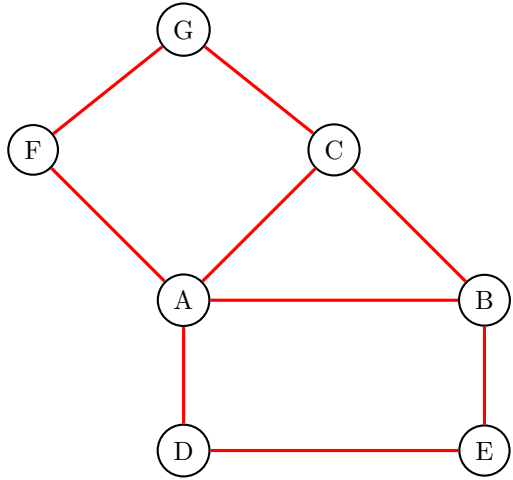
\begin{center}	
		\begin{tikzpicture}
			\begin{scope}[every node/.style={circle,thick,draw}]
				\node (A) at (0,0) {A};
				\node (B) at (4,0) {B};
				\node (C) at (2,2) {C};
				\node (D) at (0,-2) {D};
				\node (E) at (4,-2) {E};
				\node (F) at (-2,2) {F};
				\node (G) at (0,3.6) {G};
			\end{scope}		
			\begin{scope}[every edge/.style={draw=red,very thick}]
				\path (A) edge [above] node {} (B);
				\path (B) edge [above] node {} (C);
				\path (C) edge [above] node {} (A);
				\path (A) edge [above] node {} (D);
				\path (B) edge [above] node {} (E);
				\path (D) edge [above] node {} (E);
				\path (A) edge [above] node {} (F);
				\path (F) edge [above] node {} (G);
				\path (G) edge [above] node {} (C);
			\end{scope}
		\end{tikzpicture}
	\end{center}
	
	Similar analysis for $BC$, implies the existence of $H,I\notin\{A,B,C,D,E,F,G\}$, such that:
	
	
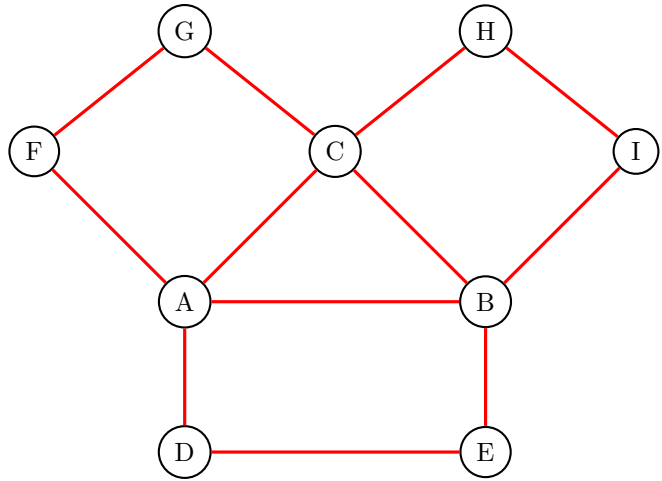
\begin{center}
		\begin{tikzpicture}
			\begin{scope}[every node/.style={circle,thick,draw}]
				\node (A) at (0,0) {A};
				\node (B) at (4,0) {B};
				\node (C) at (2,2) {C};
				\node (D) at (0,-2) {D};
				\node (E) at (4,-2) {E};
				\node (F) at (-2,2) {F};
				\node (G) at (0,3.6) {G};
				\node (H) at (4,3.6) {H};
				\node (I) at (6,2) {I};
			\end{scope}		
			\begin{scope}[every edge/.style={draw=red,very thick}]
				\path (A) edge [above] node {} (B);
				\path (B) edge [above] node {} (C);
				\path (C) edge [above] node {} (A);
				\path (A) edge [above] node {} (D);
				\path (B) edge [above] node {} (E);
				\path (D) edge [above] node {} (E);
				\path (A) edge [above] node {} (F);
				\path (F) edge [above] node {} (G);
				\path (G) edge [above] node {} (C);
				\path (C) edge [above] node {} (H);
				\path (H) edge [above] node {} (I);
				\path (I) edge [above] node {} (B);
			\end{scope}
		\end{tikzpicture}
	\end{center}

	Note that $GD$ is blue, for if otherwise, we have edge $AD$ in both quadrilaterals $\Box(DABE)$ and $\Box(DAFG)$. Similarly, $GI$ and $ID$ are blue. We have:
	
	
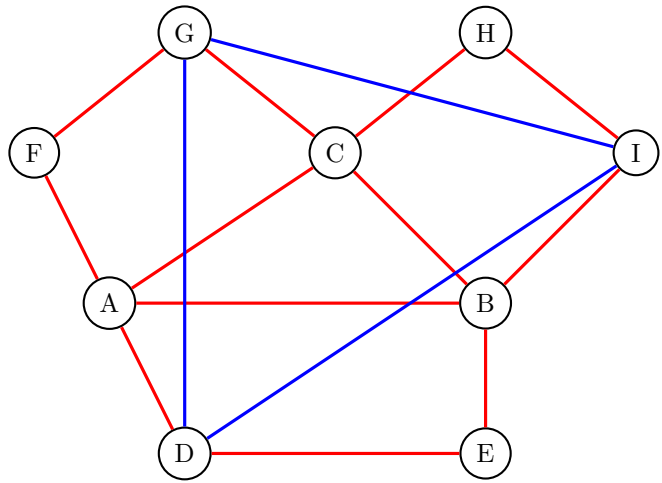
\begin{center}
		\begin{tikzpicture}
			\begin{scope}[every node/.style={circle,thick,draw}]
				\node (A) at (-1,0) {A};
				\node (B) at (4,0) {B};
				\node (C) at (2,2) {C};
				\node (D) at (0,-2) {D};
				\node (E) at (4,-2) {E};
				\node (F) at (-2,2) {F};
				\node (G) at (0,3.6) {G};
				\node (H) at (4,3.6) {H};
				\node (I) at (6,2) {I};
			\end{scope}		
			\begin{scope}[every edge/.style={draw=red,very thick}]
				\path (A) edge [above] node {} (B);
				\path (B) edge [above] node {} (C);
				\path (C) edge [above] node {} (A);
				\path (A) edge [above] node {} (D);
				\path (B) edge [above] node {} (E);
				\path (D) edge [above] node {} (E);
				\path (A) edge [above] node {} (F);
				\path (F) edge [above] node {} (G);
				\path (G) edge [above] node {} (C);
				\path (C) edge [above] node {} (H);
				\path (H) edge [above] node {} (I);
				\path (I) edge [above] node {} (B);
			\end{scope}
			\begin{scope}[every edge/.style={draw=blue,very thick}]
				\path (G) edge [above] node {} (I);
				\path (D) edge [above] node {} (G);
				\path (I) edge [above] node {} (D);
			\end{scope}
		\end{tikzpicture}
	\end{center}
	
	Now, if $AG$ is red, then $\Box(AGCB)$ is a red quadrilateral sharing $AB$ with $\Box(ABED)$, contradicting (P3). Also, if $AI$ is red, then we have a red $\Box(ACHI)$ sharing $CH$ with $\Box(BCHI)$. Therefore, $AG$ and $AI$ must be blue, and we have:
	
	
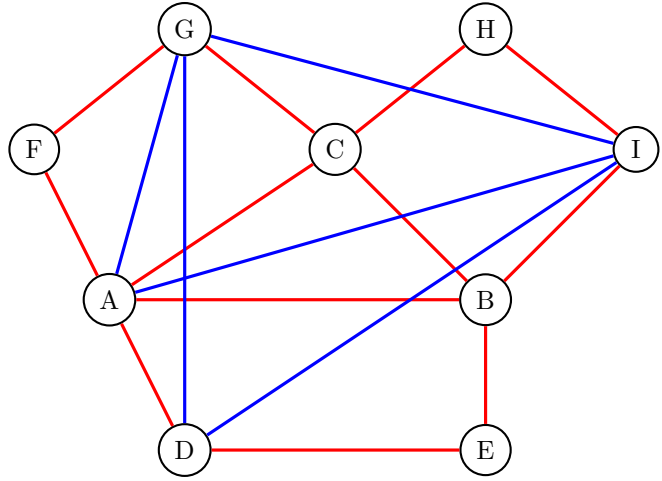
\begin{center}
		\begin{tikzpicture}
			\begin{scope}[every node/.style={circle,thick,draw}]
				\node (A) at (-1,0) {A};
				\node (B) at (4,0) {B};
				\node (C) at (2,2) {C};
				\node (D) at (0,-2) {D};
				\node (E) at (4,-2) {E};
				\node (F) at (-2,2) {F};
				\node (G) at (0,3.6) {G};
				\node (H) at (4,3.6) {H};
			\node (I) at (6,2) {I};
			\end{scope}		
			\begin{scope}[every edge/.style={draw=red,very thick}]
				\path (A) edge [above] node {} (B);
				\path (B) edge [above] node {} (C);
				\path (C) edge [above] node {} (A);
				\path (A) edge [above] node {} (D);
				\path (B) edge [above] node {} (E);
				\path (D) edge [above] node {} (E);
				\path (A) edge [above] node {} (F);
				\path (F) edge [above] node {} (G);
				\path (G) edge [above] node {} (C);
				\path (C) edge [above] node {} (H);
				\path (H) edge [above] node {} (I);
				\path (I) edge [above] node {} (B);
			\end{scope}
			\begin{scope}[every edge/.style={draw=blue,very thick}]
				\path (G) edge [above] node {} (I);
				\path (D) edge [above] node {} (G);
				\path (I) edge [above] node {} (D);
				\path (A) edge [above] node {} (G);
				\path (I) edge [above] node {} (A);
			\end{scope}
		\end{tikzpicture}
	\end{center}
	
	Implying that $IG$ is a blue edge shared by $\Delta(AGI)$ and $\Delta(DGI)$ contradicting (P2). Hence, case 1 is impossible. Q.E.D.
	
	\begin{lemma}
		Case 2 is impossible.
	\end{lemma}

	\emph{Proof}: We assume now in contradiction that case 2 holds, and w.l.o.g., that $AC$ and $BC$ lie in the same quadrilateral. Then we have:

	
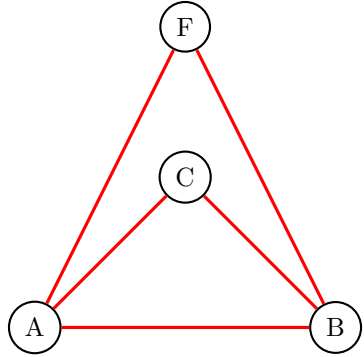
\begin{center}
		\begin{tikzpicture}
			\begin{scope}[every node/.style={circle,thick,draw}]
				\node (A) at (0,0) {A};
				\node (B) at (4,0) {B};
				\node (C) at (2,2) {C};
				\node (F) at (2,4) {F};
			\end{scope}		
			\begin{scope}[every edge/.style={draw=red,very thick}]
				\path (A) edge [above] node {} (B);
				\path (B) edge [above] node {} (C);
				\path (C) edge [above] node {} (A);
				\path (F) edge [above] node {} (A);
				\path (F) edge [above] node {} (B);
			\end{scope}
		\end{tikzpicture}
	\end{center}

	As $AB$ is part of a unique quadrilateral, there are $D,E$, such that $AB\in\Box(ABED)$. It is easy to check that $D,E\notin\{A,B,C,F\}$. We have established the fundamental structure of case 2:

	
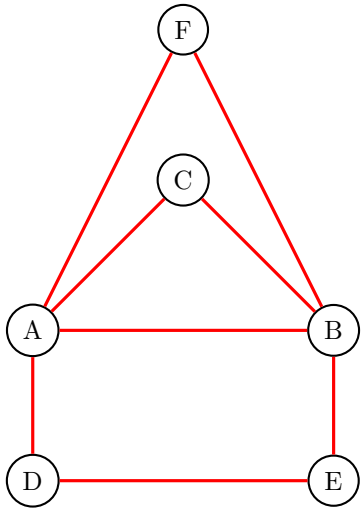
\begin{center}
		\begin{tikzpicture}
			\begin{scope}[every node/.style={circle,thick,draw}]
				\node (A) at (0,0) {A};
				\node (B) at (4,0) {B};
				\node (C) at (2,2) {C};
				\node (F) at (2,4) {F};
				\node (D) at (0,-2) {D};
				\node (E) at (4,-2) {E};
			\end{scope}		
			\begin{scope}[every edge/.style={draw=red,very thick}]
				\path (A) edge [above] node {} (B);
				\path (B) edge [above] node {} (C);
				\path (C) edge [above] node {} (A);
				\path (A) edge [above] node {} (D);
				\path (B) edge [above] node {} (E);
				\path (D) edge [above] node {} (E);
				\path (F) edge [above] node {} (A);
				\path (F) edge [above] node {} (B);
			\end{scope}
		\end{tikzpicture}
	\end{center}
		
	It is simple to see that $DC$ and $EC$ must be blue, for if $DC$ was red, then we would have quadrilaterals $\Box(DCBE)$ and $\Box(ABED)$ share $BE$. Similarly, $EC$ must be blue, then:
	
	
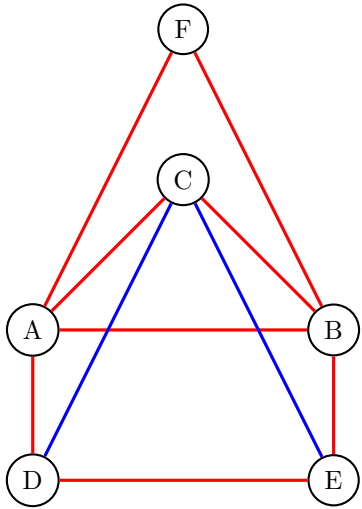
\begin{center}
		\begin{tikzpicture}
			\begin{scope}[every node/.style={circle,thick,draw}]
				\node (A) at (0,0) {A};
				\node (B) at (4,0) {B};
				\node (C) at (2,2) {C};
				\node (F) at (2,4) {F};
				\node (D) at (0,-2) {D};
				\node (E) at (4,-2) {E};
			\end{scope}		
			\begin{scope}[every edge/.style={draw=red,very thick}]
				\path (A) edge [above] node {} (B);
				\path (B) edge [above] node {} (C);
				\path (C) edge [above] node {} (A);
				\path (A) edge [above] node {} (D);
				\path (B) edge [above] node {} (E);
				\path (D) edge [above] node {} (E);
				\path (F) edge [above] node {} (A);
				\path (F) edge [above] node {} (B);
			\end{scope}
			\begin{scope}[every edge/.style={draw=blue,very thick}]
				\path (D) edge [above] node {} (C);
				\path (E) edge [above] node {} (C);
			\end{scope}
		\end{tikzpicture}
	\end{center}

	Also, we have $FD$ blue, for if otherwise, $FD$ was red, and then $\Box(FBED)$ and $\Box(ABED)$ would be red quadrilaterals sharing $BE$. Similarly $EF$ must be blue, then we have:
	
	
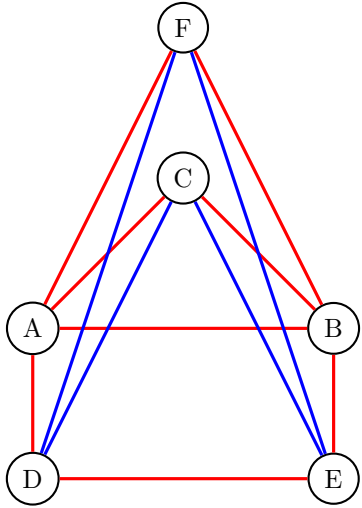
\begin{center}
		\begin{tikzpicture}
			\begin{scope}[every node/.style={circle,thick,draw}]
				\node (A) at (0,0) {A};
				\node (B) at (4,0) {B};
				\node (C) at (2,2) {C};
				\node (F) at (2,4) {F};
				\node (D) at (0,-2) {D};
				\node (E) at (4,-2) {E};
				\end{scope}		
			\begin{scope}[every edge/.style={draw=red,very thick}]
				\path (A) edge [above] node {} (B);
				\path (B) edge [above] node {} (C);
				\path (C) edge [above] node {} (A);
				\path (A) edge [above] node {} (D);
				\path (B) edge [above] node {} (E);
				\path (D) edge [above] node {} (E);
				\path (F) edge [above] node {} (A);
				\path (F) edge [above] node {} (B);
			\end{scope}
			\begin{scope}[every edge/.style={draw=blue,very thick}]
				\path (D) edge [above] node {} (C);
				\path (E) edge [above] node {} (C);
				\path (F) edge [above] node {} (D);
				\path (F) edge [above] node {} (E);
			\end{scope}
		\end{tikzpicture}
	\end{center}
	
	Thus, $FC$ must be red, for if otherwise, it would be shared by $\Delta(FDC),\Delta(FEC)$. Omitting the blue edges from the illustration for simplicity, we have:
	
	
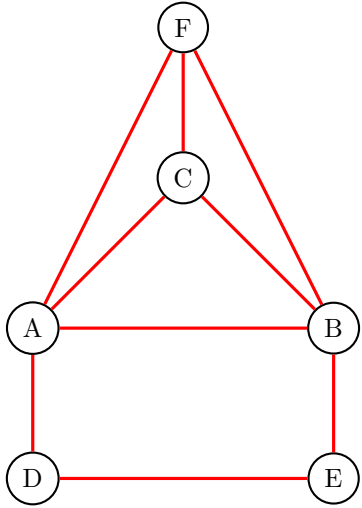
\begin{center}
		\begin{tikzpicture}
			\begin{scope}[every node/.style={circle,thick,draw}]
				\node (A) at (0,0) {A};
				\node (B) at (4,0) {B};
				\node (C) at (2,2) {C};
				\node (F) at (2,4) {F};
				\node (D) at (0,-2) {D};
				\node (E) at (4,-2) {E};
			\end{scope}		
			\begin{scope}[every edge/.style={draw=red,very thick}]
				\path (A) edge [above] node {} (B);
				\path (B) edge [above] node {} (C);
				\path (C) edge [above] node {} (A);
				\path (A) edge [above] node {} (D);
				\path (B) edge [above] node {} (E);
				\path (D) edge [above] node {} (E);
				\path (F) edge [above] node {} (A);
				\path (F) edge [above] node {} (B);
				\path (F) edge [above] node {} (C);
			\end{scope}
		\end{tikzpicture}
	\end{center}
	
	Thus, $\Box(FCBA)$ and $\Box(ABED)$ are red quadrilaterals, sharing $AB$, a contradiction to (P3). Q.E.D.

	We have thus established the impossibility of all different configurations that support a red triangle, hence there are no red triangles, proving Lemma 1.
	
\section{G has no blue quadrilaterals - proof of Lemma 2}
	We will follow similar case analysis arguments to prove the impossibility of $G$ satisfying (P2-3) having a blue quadrilateral. Assume in contradiction that $G$ satisfies (P2-3) and contains a blue quadrilateral, then we have the following substructure of $G$:
	
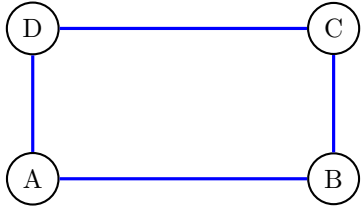
	\begin{center}
		\begin{tikzpicture}
			\begin{scope}[every node/.style={circle,thick,draw}]
				\node (A) at (0,0) {A};
				\node (B) at (4,0) {B};
				\node (D) at (0,2) {D};
				\node (C) at (4,2) {C};
			\end{scope}		
			\begin{scope}[every edge/.style={draw=blue,very thick}]
				\path (A) edge [above] node {} (B);
				\path (B) edge [above] node {} (C);
				\path (C) edge [above] node {} (D);
				\path (D) edge [above] node {} (A);
			\end{scope}
		\end{tikzpicture}
	\end{center}
	
	Can two edges, w.l.o.g., $AD$ and $AB$ reside in the same triangle? The answer is no, for if otherwise we would have:
	
	
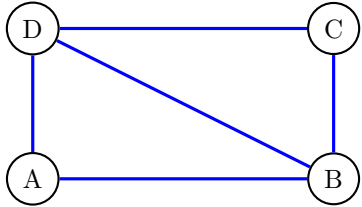
\begin{center}
		\begin{tikzpicture}
			\begin{scope}[every node/.style={circle,thick,draw}]
				\node (A) at (0,0) {A};
				\node (B) at (4,0) {B};
				\node (D) at (0,2) {D};
				\node (C) at (4,2) {C};
			\end{scope}		
			\begin{scope}[every edge/.style={draw=blue,very thick}]
				\path (A) edge [above] node {} (B);
				\path (B) edge [above] node {} (C);
				\path (C) edge [above] node {} (D);
				\path (D) edge [above] node {} (A);
				\path (D) edge [above] node {} (B);
			\end{scope}
		\end{tikzpicture}
	\end{center}
	
	Contradicting (P2), as $DB$ is shared by the triangles $\Delta(ADB)$ and $\Delta(DCB)$. Thus, we have each edge of $\Box(ABCD)$ on a different triangle. Say $DC$ is a part of $\Delta(DCE)$, for $E\notin\{A,B,C,D\}$, then we have:
	
	
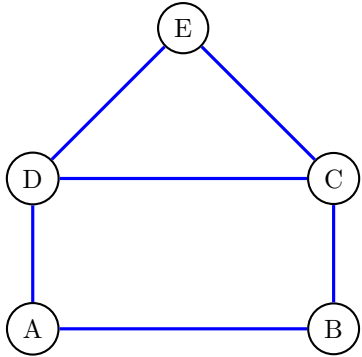
\begin{center}
		\begin{tikzpicture}
			\begin{scope}[every node/.style={circle,thick,draw}]
				\node (A) at (0,0) {A};
				\node (B) at (4,0) {B};
				\node (D) at (0,2) {D};
				\node (C) at (4,2) {C};
				\node (E) at (2,4) {E};
			\end{scope}		
			\begin{scope}[every edge/.style={draw=blue,very thick}]
				\path (A) edge [above] node {} (B);
				\path (B) edge [above] node {} (C);
				\path (C) edge [above] node {} (D);
				\path (D) edge [above] node {} (A);
				\path (E) edge [above] node {} (D);
				\path (E) edge [above] node {} (C);
			\end{scope}
		\end{tikzpicture}
	\end{center}
	
	Note that $E$ cannot be connected by a blue edge to any other vertex of $ABCD$, for if otherwise, assume w.l.o.g. that it is connected to $A$ by a blue edge, then we have $DE$ in both triangles $\Delta(ADE)$ and $\Delta(DEC)$. The previous statement, together with the fact that each blue edge belongs to a unique triangle, implies the existence of 3 other vertices, $F, G$ and $H$, forming the following structure:
	
	
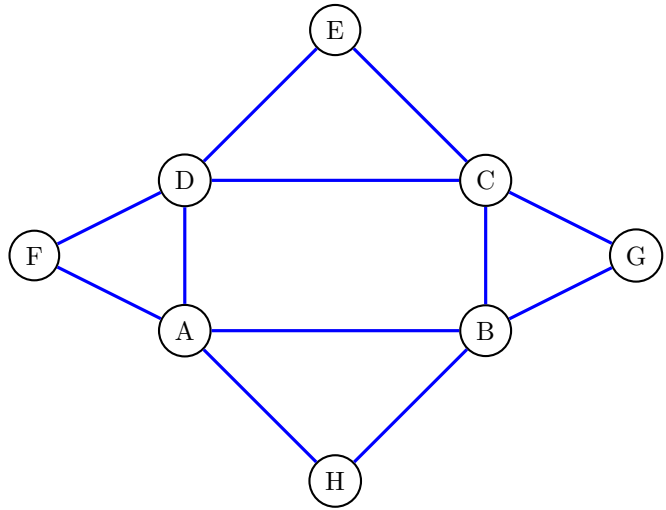
\begin{center}
		\begin{tikzpicture}
			\begin{scope}[every node/.style={circle,thick,draw}]
				\node (A) at (0,0) {A};
				\node (B) at (4,0) {B};
				\node (D) at (0,2) {D};
				\node (C) at (4,2) {C};
				\node (E) at (2,4) {E};
				\node (F) at (-2,1) {F};
				\node (G) at (6,1) {G};
				\node (H) at (2,-2) {H};
			\end{scope}		
			\begin{scope}[every edge/.style={draw=blue,very thick}]
				\path (A) edge [above] node {} (B);
				\path (B) edge [above] node {} (C);
				\path (C) edge [above] node {} (D);
				\path (D) edge [above] node {} (A);
				\path (E) edge [above] node {} (D);
				\path (E) edge [above] node {} (C);
				\path (F) edge [above] node {} (D);
				\path (F) edge [above] node {} (A);
				\path (H) edge [above] node {} (A);
				\path (H) edge [above] node {} (B);
				\path (G) edge [above] node {} (C);
				\path (G) edge [above] node {} (B);
			\end{scope}
		\end{tikzpicture}
	\end{center}
	
	To end our argument we note that $DG$ must be red, for if otherwise $CG$ would lie in both triangles $\Delta(DCG)$ and $\Delta(CGB)$. Similarly, $GH$ and $HD$ must be red, but then we have:
	
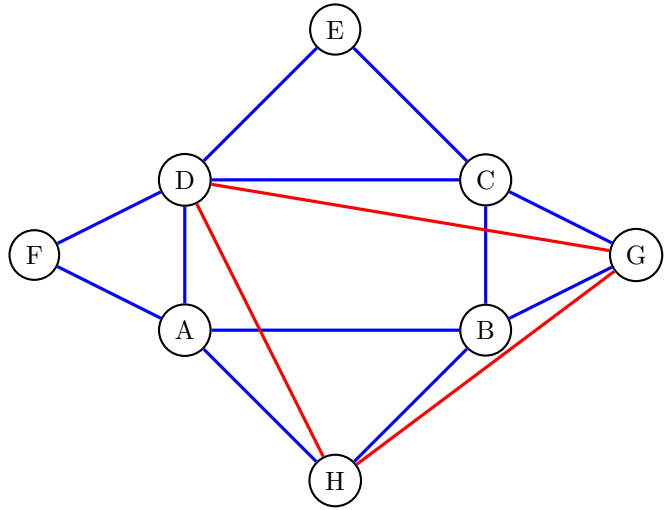
\begin{center}
		\begin{tikzpicture}
			\begin{scope}[every node/.style={circle,thick,draw}]
				\node (A) at (0,0) {A};
				\node (B) at (4,0) {B};
				\node (D) at (0,2) {D};
				\node (C) at (4,2) {C};
				\node (E) at (2,4) {E};
				\node (F) at (-2,1) {F};
				\node (G) at (6,1) {G};
				\node (H) at (2,-2) {H};
			\end{scope}		
			\begin{scope}[every edge/.style={draw=blue,very thick}]
				\path (A) edge [above] node {} (B);
				\path (B) edge [above] node {} (C);
				\path (C) edge [above] node {} (D);
				\path (D) edge [above] node {} (A);
				\path (E) edge [above] node {} (D);
				\path (E) edge [above] node {} (C);
				\path (F) edge [above] node {} (D);
				\path (F) edge [above] node {} (A);
				\path (H) edge [above] node {} (A);
				\path (H) edge [above] node {} (B);
				\path (G) edge [above] node {} (C);
				\path (G) edge [above] node {} (B);
			\end{scope}
			\begin{scope}[every edge/.style={draw=red,very thick}]
				\path (D) edge [above] node {} (G);
				\path (G) edge [above] node {} (H);
				\path (H) edge [above] node {} (D);
			\end{scope}
		\end{tikzpicture}
	\end{center}
	
	Then $\Delta(DGH)$ is a red triangle, and by Lemma 1, we reach a contradiction. Hence, there are no blue quadrilaterals.
	\newpage
\section{Proof of Theorem 2}
	Let us denote the number of blue edges by $b$ and the number of red edges by $r$. As $b$ and $r$ are respectively the edges and none edges of a simple graph $G$, we must have $r + b = {n\choose 2}$. Assuming $G$ satisfies (P2-3), by Lemma 1, the graph induced by the red edges is triangle free. It is a well known fact that a triangle free graph has at most $n^2/4$ edges. Hence, $r\le n^2/4$. Lemma 2 implies that the graph induced by the blue edges is quadrilateral free. It is known that a quadrilateral free graph has at most $(1+\sqrt{4n-3})n/4$ edges~\cite{C4}, hence $b\le\dfrac{n}{4}(1+\sqrt{4n-3})$. We have then established that:
	\begin{center}
		${n\choose 2} = r + b \le \dfrac{n^2}{4} + \dfrac{n}{4}(1+\sqrt{4n-3})$. Which implies that $n<9$. 
	\end{center}
	Q.E.D.
\section{Acknowledgments}
	We would like to thank Brendan Rooney for notifying us that we, in fact, solved the alternative interpretation of Conwel's 99-graph problem.
	
\bibliographystyle{plain}
\bibliography{documento}

\begin{thebibliography}{2}

\bibitem{CON}
John Horton Conway.
\newblock 2017 Open Problems
\url{https://oeis.org/A248380/a248380.pdf}
\newblock Accessed on 25/07/2017

\bibitem{C4}
Stasys Jukna. Extremal combinatorics: with applications in computer science.
\newblock Springer Science \& Business Media, 2011. Chapter 2, Section 2.

\end{thebibliography}
\end{document}